# Motions on n-Simplex Graphs with m-value memory

*Marc Zucker*

Department of Mathematics
Marymount Manhattan College
New York, NY 10021
mzucker@mmm.edu

Abstract: We introduce the idea of an n-simplex graph and games upon simplicial complexes. We then define moves on a labeled graph and pose the problem of whether given two labelings of a graph it is possible to change one into another via these moves. We then solve the problem for a given class of graphs. Once having found a solution for a given class of graphs we determine the number of different solutions that exist. We then use this to find an algorithm to determine whether a graph is ($n$+1)-colorable, and in particular, whether it is 3-colorable.

## 0  Introduction

The relation between graphs and games has been recognized throughout the development of Game Theory (see [3, 9, 13]). Graphs have played an integral role in the solution of coin flipping games, which are usually viewed as being upon graphs and grids. Coin flipping games, also known as **s**-games, are games in which the coins adjacent to the one being flipped are also flipped (in certain versions only the adjacent coins are flipped). Coin flipping games are also popularly known as light-switching games. That is, a game in which a switch toggles neighboring switches as well. The objective in these games is, given an initial configuration, to turn all the lights off (or coins over). **s**-games have been studied by Sutner [10, 11, 12] and by Barua and Ramakrishnan [2]. The commercially available game "Lights Out," played on a 5×5 board, has been studied by Anderson and Feil [1]. (For more on these games and their generalizations see [5, 6].) More recently, games have been looked at upon simplicial complexes [4]. Though any countable graph with vertices in n-dimensions can be viewed as a graph in two dimensions, the ability to view it in n-dimensions allows for the notions of motions and moves along the graph to have a more intuitive sense. The use of simplicial complexes is then a natural generalization.

In this paper we introduce a game in which each move, called a '*push*', changes the labeling of every vertex in a clique (not necessarily maximal). Since non-maximal cliques can be viewed as n-simplexes, we opt for the more visual notion of an n-simplex. Due to the label changing action of a push there is some similarity to **s**-games; however, in many senses the similarity goes no further than the notion of a flipping game. In Sec. 4 we give a criterion to determine whether given two boards, which we define and call n-simplex graphs, it is possible to change one board into the other by a series of these moves. In Sec. 5 we determine (given that a solutions exists) how many different solutions there are. We use these results to find an algorithm, in Sec. 7, to determine a sufficient criterion for whether a graph is $n+1$-colorable, and in specific, whether a planar graph is 3-colorable.

# 1   Definitions

Following Hatcher [7], we first define an n-simplex and a simplicial complex topologically and then indicate their natural corresponding graphical interpretation.

**Definition 1.1** *Given any set* $V = \{v_0, v_1, ..., v_n\}$ *of* $n+1$ *points in* $\mathbf{R}^N$, *such that the differences* $v_1 - v_0, v_2 - v_0, ..., v_n - v_0$ *are linearly independent, the* n-simplex *with vertices* V *is the convex hull of V, i.e. the set of all points of the form* $t_0 v_0 + t_1 v_1 + ... + t_n v_n$, *where* $\sum_{i=0}^{n} t_i = 1$ *and* $t_i \geq 0$ *for all i.*

**Definition 1.2** *A simplicial complex* $\Delta$ *on a finite set* V *is a collection of subsets of V such that*

*i)* $\{v\} \in \Delta$ *for all* $v \in V$.
*ii) if* $F \in \Delta$ *and* $G \subseteq F$ *then* $G \in \Delta$.

*The members of* $\Delta$ *are called* simplices *or* faces, *and the elements of V are called vertices.*

The use (and advantage) of topological notions is that it allows us to have n-simplex graphs upon any compact connected $q$-manifold, with $q \geq n$.

Viewed purely from a graph theoretic perspective, an n-simplex is simply a complete graph, $U_i$, such that $|V(U_i)| = n+1$. Thus the set of vertices, $V(U_i)$, uniquely determines the n-simplex. And a simplicial complex is merely a graph $G = \bigcup U_j$ with each $U_j$ complete. Due to the natural mapping from the topological definition of simplicial complexes to the graph-theoretic definition, we use them interchangeably, adopting the notion of an n-dimensional graph, $G$, comprised of n-simplexes. It should be noted that the connection to the topological notion has certain limitations since we are not concerned with the topology of n-simplexes or with viewing them as convex hulls, but rather with their vertex sets, edge sets, and faces. Therefore, except for the intuition

that it brings with it, the ideas are primarily graph-theoretic, and should be viewed as such.

We now define the notion of an n-simplex graph based upon the notions of a simplicial complex.

**Definition 1.3** *G is said to be an* n-simplex graph *if:*

    i)    *G is a simplicial complex.*
    ii)    *For any $F \in G$ with $|F| < n+1$ there exists a $K \in G$ such that $F \subset K$.*

Alternately an n-simplex *graph* can be defined (and viewed) as a graph $G = \bigcup U_j$ where each $U_j$ is complete and $|V(U_j)| = n+1$.

We take it as a general assumption throughout this paper that all graphs are of finite size.

## 2    *Region-paths and region-connected graphs*

The building blocks of n-simplex graphs are the n-simplexes. Topologically they are the convex hulls of the set of vertices. It is therefore natural that we shall refer to these n-simplexes as the regions of *G*.

We now extend the familiar notions of adjacent vertices and paths along vertices to that of regions. Adjacent regions are defined in the same vein as that of adjacent vertices.

**Definition 2.1** *Two n-simplexes, $S_i$ and $S_j$ are said to be* adjacent *if $S_i \cap S_j$ is an (n-1)-simplex.*

**Definition 2.2** *A* region-path *from $S_i$ to $S_j$ is a set of n-simplexes $\{S_i, S_{i+1},...,S_j\}$ such that $S_k$ is adjacent to $S_{k+1}$ for $i \leq k < j$.*

Let $L(G): V(G) \to Z_n$ be a labeling of the vertices of a graph $G$ from the set $\{0, 1,...,n\}$. We call a move a *push* if it is a function $f_{S_i}: L_1(G) \to L_2(G)$, acting on an n-simplex $S_i = \{v_0, v_1,...,v_n\}$, such that $f_{S_i}[L_1(v_j)] = L_1(v_j) + 1 \pmod{n}$ for $v_j \in S$.

**Definition 2.3** *A graph is said to be* region-connected *if given any two n-simplexes* $S_1, S_2 \in G$ *there exists a region-path connecting them.*

Pictorially, a *region-connected* n-simplex graph can be viewed as a graph formed by pasting n-simplexes together by their $(n-1)$-simplexes. That is, the intersection of two n-simplexes which are pasted together is an $(n-1)$-simplex. For example, we could form a region-connected 2-simplex graph by pasting triangles together by their edges, or form a region-connected 3-simplex graph by pasting tetrahedrons together by their triangles.

## 3 An invariant under pushes

We pose the following two general problems concerning motions, or re-labelings, of graphs:

**Question 3.1** *Let G be any region-connected n-simplex graph. Given two labelings, $L_1(G)$ & $L_2(G)$, does there exist a series of pushes with which we can change $L_1(G)$ into $L_2(G)$?*

**Question 3.2** *Given that a solution to Question 3.1 exists, how many different solutions are there?*

We leave the solution to Question 3.2 until Sec. 5.

We wish to find a class of n-simplex graphs under which our question is solvable. We claim that if $G$ is a region-connected n-simplex graph then the question is solvable if $c(G) = n+1$ (not that it is affirmative).

We start by finding a value of $L(G)$ which is invariant under pushes.

Let $Z_m$ be the labeling set for $L_1$ & $L_2$.

Since $c(G) = n+1$, $G$ can be colored with the set $\{i_0, i_1, ..., i_n\}$, where $i_k$ is as follows:

$$i_k = \begin{bmatrix} I_k & & & 0 & 0 \\ & \ddots & & & 0 \\ & & e^{\frac{i2p}{m}} & & \\ 0 & & & \ddots & \\ 0 & 0 & & & I_{n-k-1} \end{bmatrix} \quad 0 \leq k < n, \quad i_n = \begin{bmatrix} e^{\frac{i2(m-1)p}{m}} & & 0 \\ & \ddots & \\ 0 & & e^{\frac{i2(m-1)p}{m}} \end{bmatrix} = e^{\frac{i2(m-1)p}{m}} I_n$$

We have $i_0^m = i_1^m = \cdots = i_{n-1}^m = i_n^m = i_0 i_1 \cdots i_n = I$.

Given that $\{i_0, i_1, \ldots, i_{n-1}\}$ is linearly independent and that the order of each $i_j$ is $m$, we find that $\langle i_0, i_1, \ldots, i_{n-1} \rangle \approx \underbrace{Z_m \times Z_m \times \ldots \times Z_m}_{n-times}$.

Let $i(v_j) \in \{i_0, i_1, \cdots, i_n\}$ be the coloring of the vertex $v_j \in V(G)$ and $l(v_j) \in L(G)$ its label. Assign the value $i^{l(v_j)}(v_j)$ to the vertex $v_j$, for each $v_j \in V(G)$.

Let $P[L(G)] = \prod_{v_j \in V(G)} i^{l(v_j)}(v_j)$.

**Lemma 3.2**    $P[L(G)]$ *is invariant under pushes.*

**Proof**    Let $\{i_0^{a_0}, i_1^{a_1}, \cdots, i_n^{a_n}\}$ be the set of values assigned to the vertices of an arbitrary n-simplex in $G$. A push on this n-simplex would send $i_j^{a_j} \mapsto i_j^{a_j+1}$ $\forall j \in \{0,1,\ldots,n\}$. Thus we have:

$$i_0^{a_0} i_1^{a_1} \cdots i_n^{a_n} \to i_0^{a_0+1} i_1^{a_1+1} \cdots i_n^{a_n+1} = i_0^{a_0} i_1^{a_1} \cdots i_n^{a_n} (i_0 i_1 \cdots i_n) = i_0^{a_0} i_1^{a_1} \cdots i_n^{a_n} \quad \text{(since } i_0 i_1 \cdots i_n = I \text{)}.$$

(Since the push acts only upon a given n-simplex, the remainder of $G$ remains unchanged.)

$P[L(G)]$ is therefore invariant under pushes.

## 4    Solutions for a specific class of graphs

We now show a class of graphs under which our question is solvable and present an algorithm for finding a series of pushes when a solution exists.

**Theorem 4.1** *Let G be a region-connected n-simplex graph with $c(G) = n+1$. Then there exists a set F, of pushes, such that $F[L_1(G)] = L_2(G)$ iff $P[L_1(G)] = P[L_2(G)]$.*

**Proof** By Lemma 3.2, $P[L(G)]$ is invariant under pushes; the necessity of the equality therefore follows. (The condition of region-connectedness is not actually required for this direction.)

To show the sufficiency of the equality we will construct a solution. Since $G$ is region-connected, we can find n+1 region-paths in G such that:

(i) $s_{1,i_j}, s_{2,i_j}, \cdots, s_{t,i_j}$ is a region-path of n-simplexes with $s_{(2k-1),i_j} = \{v_0, v_1, \cdots, v_n\}$, $s_{(2k),i_j} = \{v_1, \cdots, v_n, v_{n+1}\}$, then $s_{(2k-1),i_j} \cap s_{(2k),i_j} = \{v_1, \cdots, v_n\}$, $k = \{1, 2, \ldots\}$, where $i_j$ is the coloring of both $v_0$ and $v_{n+1}$.

(ii) Given any $v \in V(G)$ with coloring $i_j$, $v$ is in $s_{r,i_j}$ for some $r$.

(iii) $s_{t,i_0} = s_{t,i_1} = \cdots = s_{t,i_n}$. I.e. all paths end with the same n-simplex.

(No other conditions exist for these paths; thus it is possible for a region-path to cross or retrace itself.)

Let us now assume that two labelings, $L_1$ & $L_2$, differ in only one n-simplex. Let $f_{1,i_0}^{k_{1,i_0}}$ be a push to such a power that the vertex in this n-simplex colored $i_0$ attains the same labeling as in $L_2(G)$. We claim that $L_1(G) = L_2(G)$. Since $P[L_1(G)] = P[L_2(G)]$ (by assumption) we have $i_0^{a_{1,0}} i_1^{a_{1,1}} \cdots i_n^{a_{1,n}} = i_0^{a_{2,0}} i_1^{a_{2,1}} \cdots i_n^{a_{2,n}}$. And given that $i_0^{a_{1,0}} = i_0^{a_{2,0}}$, we have $A_{ij} = i_1^{a_{1,1}}, \cdots, i_n^{a_{1,n}} = i_1^{a_{2,1}}, \cdots, i_n^{a_{2,n}} = B_{ij}$. As well, since $a_{1,1} = e^{\frac{i2a_{1,n}(m-1)p}{m}} = e^{\frac{i2a_{2,n}(m-1)p}{m}} = b_{1,1}$, we have $\mathbf{a}_{1,n} = \mathbf{a}_{2,n}$. That the remaining labels are also equal can be seen as a result of their linear independence. We therefore have that $\mathbf{a}_{1,j} = \mathbf{a}_{2,j} \ \forall j \in \{0,1,\cdots,n\}$.

Let us now assume that the labelings differ arbitrarily. We form $n+1$ sequences of pushes as follows: For the coloring $i_j$ we have the sequence $f_{1,i_j}^{p_{1,i_j}} f_{2,i_j}^{m-p_{1,i_j}} f_{3,i_j}^{p_{3,i_j}} f_{4,i_j}^{m-p_{3,i_j}} \cdots$, where the pushes act upon the n-simplexes of the region-paths given above, and where the $(2k-1)^{st}$ term is of the form $f_{(2k-1),i_j}^{p_{(2k-1),i_j}}$ and the $(2k)^{th}$ term is of the form $f_{(2k),i_j}^{m-p_{1,i_j}}$, for $k = \{1,2,\ldots\}$, and where $p_{(2k-1),i_j}$ is the power necessary so that if $v$ is the vertex in the n-simplex colored $i_j$, then $l_1\left[f_{(2k-1),i_j}^{p_{(2k-1),i_j}}(v)\right] = l_2(v)$. I.e. the labeling of that vertex is the same as in $L_2(G)$.

The last term in each of the sequences is either of the form $f_{(2k-1),i_j}^{P_{(2k-1),i_j}}$ or $f_{(2k),i_j}^{m-P_{(2k-1),i_j}}$, depending upon whether there is an odd or an even number of elements in the region-path.

Note that since $s_{(2k-1),i_j} \cap s_{(2k),i_j} = \{v_1, \cdots, v_n\}$, $f_{(2k-1),i_j}^{P_{(2k-1),i_j}} f_{(2k),i_j}^{m-P_{(2k-1),i_j}}$ raises the values of all those vertices not colored $i_j$ by $m$. Therefore, all vertices, except those colored $i_j$, remain unchanged.

Each sequence will therefore change the labelings of all those vertices of $G$ colored $i_j$, except, possibly, for that of the last n-simplex of the path (upon which the sequence acts). This being true for each $i_j$, we need only concern ourselves now with this last n-simplex, which, by construction, is the same for each path (the rest of $G$ can therefore be ignored). Having already proven the theorem for that case our proof is complete.

## 5 Board classes and solution sizes

Our objective in this section is to show the number of different solutions that exist when there is a solution. To accomplish this we show that the set of labelings of $G$ form equivalence classes. This we do because we wish to demonstrate that the sizes of all these equivalence classes are the same.

Let us add to our collection of pushes acting on $G$, $n$ additional elements acting only upon one vertex out of the $n+1$ vertices in some given n-simplex. Without loss of generality we can assume that these new elements all act upon the same n-simplex. Given $f = \begin{bmatrix} e^{\frac{i2\mathbf{p}}{m}} & & 0 \\ & \ddots & \\ 0 & & e^{\frac{i2\mathbf{p}}{m}} \end{bmatrix}$ (i.e. multiplication of the values of the vertices in an n-simplex by this matrix gives the result of a push on that n-simplex), we let

$h_k = \begin{bmatrix} I_k & & 0 & 0 \\ & \ddots & & 0 \\ & & e^{\frac{i2\mathbf{p}}{m}} & \\ 0 & & \ddots & \\ 0 & 0 & & I_{n-k-1} \end{bmatrix}$  $i = \{0, 1, ..., n-1\}$, be our $n$ new moves. Since $P[L(G)]$ can take on any value (due to the addition of these new moves), it is now possible to change from any given labeling of $G$ to any other labeling of $G$.

We now form classes of labelings of $G$ such that two labelings, $K$ and $L$, are in the same class if $P[K(G)] = P[L(G)]$.

**Definition 5.1** Let $L_1(G)$ & $L_2(G)$ be two labelings of a region-connected n-simplex graph, G, with the condition that $c(G) = n+1$. Then we say that $L_1(G)$ is label-equivalent to $L_2(G)$, written as $L_1(G) \sim L_2(G)$, if $P[L_1(G)] = P[L_2(G)]$.

**Lemma 5.2** The relation $L_1(G) \sim L_2(G)$ is an equivalence relation.

**Proof** Using Theorem 4.1, it can be shown that the conditions for an equivalence relation are easily satisfied.

Since an equivalence relation creates a decomposition of the set into mutually disjoint subsets (see [8] for a simple proof), we have the following:

**Corollary 5.3** The equivalence relation ~ provides a decomposition of the set of all labelings of G into distinct (mutually disjoint) equivalence classes.

Let $m$ be the size of the labeling set of $G$ and $n$ be such that $c(G) = n+1$. We now have the further fact, that

**Corollary 5.4** There are $m^n$ distinct equivalence classes.

**Proof** The proof follows easily from Theorem 4.1 and the fact that $P[L(G)]$ can be anyone of $m^n$ different values.

Let $h_i^{a_i}$ be a mapping from one class of labelings to another class, then if $K$ and $L$ are two labelings in some class, $h_i^{a_i} K$ and $h_i^{a_i} L$ will be two labelings in this other class. However, it is easily shown that $h_i^{a_i} K = h_i^{a_i} L$ if and only if $K = L$. Therefore there is a one to one relationship between the different equivalence classes, and we have

**Corollary 5.5** The equivalence relation ~ divides the set of labeled graphs into equivalence classes of equal size.

By Corollary 5.4 we know that there are $m^n$ different values that $P$ can take. Letting $|V(G)| = v$, there are $m^v$ different labelings of $G$. Since, by Corollary 5.5, the sizes of the equivalence classes are the same, there are $m^{v-n}$ different labelings for each class. Therefore, we have that

**Corollary 5.6** *There are exactly $m^{v-n}$ elements in each equivalence class formed by the equivalence relation ~.*

We now form classes of pushes that all act in a similar manner on a labeling. Let $R(G)$ be the set of all n-simplexes (i.e. regions) in $G$.

**Definition 5.7** *Let f and g be two words in $R(G)$ and L some labeling of G, then we say that f is* congruent *to g, written as $f \circ g$, if $f[L(G)] = g[L(G)]$.*

**Lemma 5.8** *The relation $f \circ g$ is an equivalence relation.*

**Proof** Follows from the definition of an equivalence relation.

As before, since an equivalence relation creates a decomposition of the set into mutually disjoint subsets (see [8]), we have the following:

**Corollary 5.9** *The equivalence relation $\circ$ provides a decomposition of the set of all words in $R(G)$ into distinct (mutually disjoint) equivalence classes.*

**Corollary 5.10** *The equivalence relation $\circ$ divides $R(G)$ into equivalence classes of equal size.*

**Proof** We form a mapping $j : K \to L$ from one equivalence class, $K$, into another equivalence class, $L$, by $f \mapsto mf$, where $m$ is such that $mf \in L$. It can easily be shown that $mf \equiv mg$ iff $f \equiv g$, and $mf = mg$ iff $f = g$. Therefore, since $K$ and $L$ were arbitrary, there is a one to one correspondence between equivalence classes.

Let $|R(G)| = r$, then there are $m^r$ different sets of moves possible on $G$. Since we know by Corollary 5.10 that the classes of moves are each the same size, and by Corollary 5.6 that there are $m^{v-n}$ different classes of labelings upon which these moves act, we have that there are $\frac{m^r}{m^{v-n}} = m^{r-v+n}$ different sets of pushes for each class of labelings. Thus we have proved

**Theorem 5.10** *Given a graph G, as in Theorem 4.1, the number of solutions that exist with which one labeling can be changed into another labeling, so long as a solution exists, is $m^{r-v+n}$.*

# 6 Examples of games on 2-simplexes

**Example 6.1** As an example, imagine a board of coins (or disks with each side a different color), which are tightly packed. That is, coins laid out so that the board is made up of little triangles of coins (i.e. each triangle consisting of three coins). A general board of this type can be formed by first taking three coins, each touching the other two, then adding new coins, one at a time, so that each new coin touches at least two other coins (which each touch the other). For example, the board can be shaped hexagonally or triangularly (each row having one more coin than the previous row). Our question would now be stated as follows: Given any two of these boards (and an arrangement of the coins in heads and tails for each), does there exist a set of pushes such that one board can be changed into the other?

Since $c(G) = 3$ (which can be easily demonstrated) a solution is now easily determined (with the existence of a solution depending on whether or not $P[L_1(G)] = P[L_2(G)]$).

For example, say we would like to change 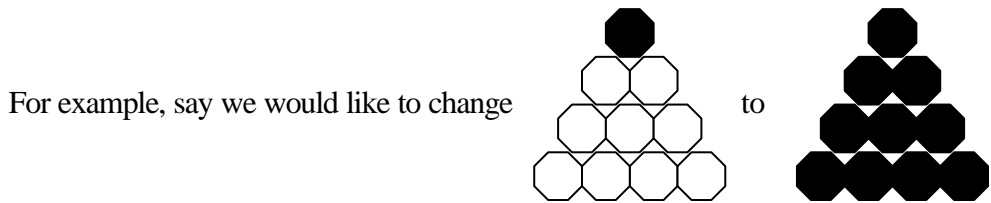 to

Since the values are equal (i.e. $P[L_1(G)] = P[L_2(G)]$, if we would 3-color it), there exists a solution, and since $m^{r-v+n} = 2^{9-10+2} = 2$, we find that there are in fact two. Upon investigation these solutions are found to be:

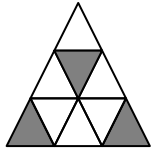 and 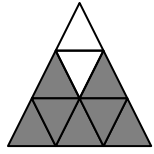

Where we obtain 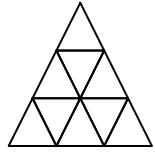 by connecting the disks as such 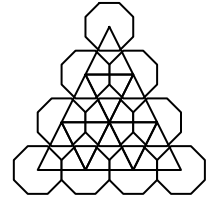

(I.e. we associate the center of each disk with a vertex and form edges between these vertices. The vertices of the colored triangles are then flipped, that is, the pushes act upon the colored triangles).

**Example 6.2** Giving the boards a memory value, we can develop Example 6.1 further. That is, we can assign a number such that the coins turn from head to tails (or tails to heads) only after a push has acted upon a vertex that number of times; for example 3 times (i.e. when a vertex is labeled 0, 1, or 2 the coin will be heads, and when it is labeled 3, 4, or 5 it will be tails). Then, since $c(G) = 3$, we can change one board into the other, as long as, again, $P[L_1(G)] = P[L_2(G)]$.

The boards used for these games can also be structured upon manifolds having a genus other than 0. For example, imagine a 2-simplex graph upon a torus (or Klein bottle), with pushes acting only upon those 2-simplexes whose convex hull (i.e. a 2-simplex in the topological sense of the term) is simply connected. This extends naturally to any compact connected n-manifold. We can therefore imagine (trivially) a 1-simplex graph on a line in 1-dimension; a simple closed curve in 2-dimensions; a knot in 3-dimensions; etc. (And similarly for other n-simplex graphs.)

These games can also be extended to two (or more) person games in a manner such as the following. Given a triangular board made up of heads, as described above (or similarly with a hexagonal board), two players are each assigned a corner. The object of the game is then to form a path of tails to the third corner, with the players alternating pushes and a win accruing for whomever succeeds first.

# 7   An algorithm for (*n*+1)-colorability

Since if $c(G) = n+1$ there are $m^n$ different classes (by Corollary 5.4) (given that $G$ is a region-connected n-simplex graph), as we move through all $m^r$ different moves possible, we obtain $\frac{1}{m^n}$ of the possible labelings for *G*. If, however, the number of different classes is $m^{n-1}$, then in the worst case scenario, it is possible that in the first $\frac{1}{m^n}$ of possible moves every labeling of one of the classes is obtained, and in the next $\frac{1}{m^n}$ of possible moves every labeling of a different class is obtained. In which case it would be possible to determine whether or not $G$ is $(n+1)$-colorable in at most $\frac{m^r}{m^n}+1 = m^{r-n}+1$ moves. (The situation gets only better if the number of different classes for $G$ is less than $m^{n-1}$.)

Therefore, we have only to show that if $G$ is not $(n+1)$-colorable then the number of classes of labelings of $G$ decreases by a factor of $m$.

Let $G$ be a graph that is more than $(n+1)$-colorable. Then there is a proper sub-graph of $G$, say $G'$, such that it is maximally $(n+1)$-colorable, with $V(G') = V(G)$. (That is, if we would add any edge from $E(G) - E(G')$ to $E(G')$, $G'$ would no longer be $(n+1)$-colorable.) Therefore, every edge in $E(G) - E(G')$ must connect two vertices of the same color. What results is that the value of the graph (i.e. $P[L(G')]$) would change if acted upon by a push; and the labeling would therefore no longer be of the same class. The actual value would change by $i_k \cdot i_{k+1}^{m-1}$, where $i_k$ is the color of both ends of the new edge. (Without loss of generality we can allow $i_{k+1}$ be the color that is missing.) This is so since a push in this case would result in a change of $i_0 \cdot i_1 \cdots i_k \cdot i_k \cdot i_{k+2} \cdots i_{n-1}$. But $i_0 \cdot i_1 \cdots i_k \cdot i_{k+2} \cdots i_{n-1} = i_0 \cdot i_1 \cdots i_k \cdot i_{k+1} \cdot i_{k+1}^{m-1} \cdot i_{k+2} \cdots i_{n-1} = (i_0 \cdot i_1 \cdots i_{n-1})i_{k+1}^{m-1} = i_{k+1}^{m-1}$, therefore $i_0 \cdot i_1 \cdots i_k \cdot i_k \cdot i_{k+2} \cdots i_{n-1} = i_k \cdot i_{k+1}^{m-1}$. However, $i_k \cdot i_{k+1}^{m-1}$ is of order $m$ (where $m$ is the size of the labeling set), so this new edge introduces $m$ new values. Thus the number of classes of labels is decreased by a factor of $m$ (since every class of labelings are now associated to $m$ other classes). For every new edge now added the number of classes could decrease by a factor of $m$, depending upon whether or not each new edge creates a new value when a push is applied to its n-simplex. An $(n+2)$-colorable graph might therefore have $m^{n-1}, m^{n-2}, \ldots, 1$ labeling classes. We have shown that

**Lemma 7.1**  *An $(n+2)$-colorable (region-connected n-simplex) graph has at most $m^{n-1}$ labeling classes.*

For example, the following 4-colorable 2-simplex graph has two labeling classes. 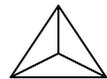

While this 4-colorable 2-simplex graph has only one labeling class. 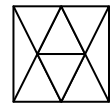

(I.e. it is possible to achieve any labeling of this graph with only pushes.)

Given a graph $G$, we form a new graph, $G'$, with $V(G') = V(G)$ and $E(G') = E(G) + \{$enough edges so that $G'$ would be a region-connected n-simplex graph$\}$. Then

**Corollary 7.2** $m^{r-n} + 1$ *moves is sufficient to determine whether a graph G is (n+1)-colorable, where* $r = R(G')$. *If, in addition, G is a region-connected n-simplex graph, then it is necessary as well.*

Since $m$ is independent of $G$ (i.e. dependent upon the labeling set only), it can be chosen arbitrarily. We can therefore improve our bound simply by letting $m = 2$.

Given an n-simplex graph $G = \bigcup_{i=1}^{N} C_i$ with each $C_i$ region-connected and $C_i \cap C_j$, $i \neq j$, a subset of an n-simplex, we form a new graph $G'$ by associating to each $C_i \subset G$ a vertex $c_i \subset V(G')$, where $c_i c_j \in E(G')$ given that $C_i \cap C_j$, $i \neq j$, is a subset of an n-simplex.

**Corollary 7.3** *Let* $G = \bigcup_{i=1}^{N} C_i$ *be an n-simplex graph where each* $C_i$ *is region-connected and where* $C_i \cap C_j$, $i \neq j$, *is a subset of an n-simplex. Then the condition in Corollary 7.2 is both necessary and sufficient for determining whether or not G is (n+1)-colorable if it's associated graph,* $G'$, *has no cycles.*

**Proof** If $G$ is $(n+1)$-colorable, then given any $C_i$ and $C_j$ as stated, since their intersection is a subset of an n-simplex, the coloring of $C_j$ can be chosen based upon that of $C_i$ (owing to the $(n+1)$-colorability of $C_j$). (The fact that there are no cycles in $G'$ guarantees this ability.) Thus it is possible to choose the coloring of $C_j$ (even if $G$ is planar), such that we can add enough edges into $E(G)$, whose vertices are colored differently, so that $G$ now becomes region-connected without changing the $(n+1)$-colorability of $G$.

While the bounds are given for arbitrary graphs, given specific conditions we can improve upon this bound. For example, for planar graphs, since $r \leq 2v - 4$, we have as a sufficient condition for 3-colorability that we need only try $2^{2v-7} + 1$ moves.

# Acknowledgements

I would like to give special thanks to Michael Anshel for all of his guidance, suggestions, comments, and critiques and to Joseph Malkevitch for his helpful conversations and suggestions. I would also like to thank Steven Wat for all of his helpful suggestions and conversations as well as the faculty of the Graduate Center of The City University of New York.